\documentclass[]{article}
\usepackage{amssymb,enumerate}
\usepackage{amsmath}
\usepackage{graphics}
\usepackage[utf8]{inputenc}

\newtheorem{theorem}{Theorem}[section]
\newtheorem{corollary}[theorem]{Corollary}

\newtheorem{remark}[theorem]{Remark}

\newtheorem{conjecture}[theorem]{Conjecture}
\title{Bounds for some entropies and special functions}
\author{Adina  B\u{a}rar\thanks{Technical University of Cluj-Napoca, Department of Mathematics, Memorandumului Street 28, 400114, Cluj-Napoca, Romania}, Gabriela Raluca Mocanu\thanks{Romanian Academy, Cluj-Napoca Branch, Astronomical Institute, Cire\c{s}ilor Street 19, 400487, Cluj-Napoca, Romania,
gabriela.mocanu$@$academia-cj.ro}, Ioan Ra\c{s}a\thanks{Technical University of Cluj-Napoca, Department of Mathematics \newline, Memorandumului Street 28, 400114, Cluj-Napoca, Romania}}

\date{}
\begin{document}

\maketitle

Subject Class: 94A17, 33E30, 33C05, 33C45.
\newline \newline
\textbf{Abstract}
\newline
We consider a family of probability distributions depending on a real parameter and including the binomial, Poisson and negative binomial distributions. The corresponding index of coincidence satisfies a Heun differential equation and is a logarithmically convex function. Combining these facts we get bounds for the index of coincidence, and consequently for R\'{e}nyi and Tsallis entropies of order $2$.

\section{Introduction}

For $c\in \mathbb{R}$, let $I_c:=\left [ 0, -\frac{1}{c} \right ]$ if $c<0$, and $I_c := [0,+\infty )$ if $c \geq 0$.

Let $a\in \mathbb{R}$ and $k\in \mathbb{N}_0$; the binomial coefficients are defined as usual by
\begin{equation*}
{a\choose k}:= \frac{a(a-1)\dots (a-k+1)}{k!}
\end{equation*}
if $k\in \mathbb{N}$, and ${a\choose 0}:=1$.

Consider also a real number $n>0$ such that $n>c$ if $c\geq 0$, and $n=-cl$ for some $l\in \mathbb{N}$ if $c<0$.

For $k\in \mathbb{N}_0$ and $x\in I_c$ define
\begin{equation*}
p_{n,k}^{[c]}(x):={-\frac{n}{c}\choose k}(-cx)^k (1+cx)^{-\frac{n}{c}-k}, \mbox{ if } c\neq 0,
\end{equation*}
\begin{equation*}
p_{n,k}^{[0]}(x):= \lim _{c\to 0}p_{n,k}^{[c]}(x) = \frac{(nx)^k}{k!}e^{-nx}.
\end{equation*}

These functions were intensively used in Approximation Theory: see \cite{3}, \cite{5}, \cite{16} and the references therein.

In particular,
\begin{equation*}
\sum _{k=0}^\infty p_{n,k}^{[c]}(x)=1,
\end{equation*}
so that $\left ( p_{n,k}^{[c]}(x) \right )_{k\geq 0}$ is a parameterized probability distribution.

Its index of coincidence (see \cite{4}) is
\begin{equation*}
S_{n,c}(x):= \sum _{k=0}^\infty \left ( p_{n,k}^{[c]}(x)\right) ^2, \quad x\in I_c.
\end{equation*}

The R\'{e}nyi entropy of order 2 and the Tsallis entropy of order 2 are given, respectively, by (see \cite{13}, \cite{15})
\begin{equation*}
R_{n,c}(x):=-\log S_{n,c}(x); \quad T_{n,c}(x):=1-S_{n,c}(x),
\end{equation*}
while the associated Shannon entropy is
\begin{equation*}
H_{n,c}(x):=-\sum _{k=0}^\infty p_{n,k}^{[c]}(x) \log p_{n,k}^{[c]}(x), \quad x\in I_c.
\end{equation*}

The cases $c=-1$, $c=0$, $c=1$ correspond, respectively, to the binomial, Poisson, and negative binomial distributions; see also \cite{9}, \cite{10}.

It was proved in \cite{8}, \cite{11} that the index of coincidence $S_{n,c}$ satisfies the Heun differential equation
\begin{eqnarray}
x(1+cx)(1+2cx)S_{n,c}''(x)+\left ( 4(n+c)x(1+cx)+1 \right)S_{n,c}'(x) +\nonumber \\ + 2n(1+2cx)S_{n,c}(x)=0, \quad x\in I_c. \label{eq:1.1}
\end{eqnarray}

It was conjectured in \cite{7} and proved in several papers (for details see \cite{1}, \cite{d},

\cite{6}, \cite{8}, \cite{11}, \cite{12} and the references given there) that $S_{n,c}$ is a convex function, i.e.,
\begin{equation}
S_{n,c}''(x)\geq 0, \quad x\in I_c. \label{eq:1.2}
\end{equation}

It is easy to combine \eqref{eq:1.1} and \eqref{eq:1.2} in order to get
\begin{equation}
S_{n,c}(x) \leq \left ( 4(n+c)x(1+cx)+1 \right )^{-\frac{n}{2(n+c)}}, \quad x\in I_c. \label{eq:1.3}
\end{equation}

Particular cases and related results can be found in \cite{8}, \cite{11}. Let us remark also that $S_{n,c}(0)=1$.

The upper bound for $S_{n,c}$, given by \eqref{eq:1.3}, leads obviously to lower bounds for the R\'{e}nyi entropy $R_{n,c}$ and the Tsallis entropy $T_{n,c}$.

The following conjecture was formulated in \cite{8} and \cite{11}:

\begin{conjecture}\label{conj:1.1} For $c\in \mathbb{R}$, $S_{n,c}$ is a logarithmically convex function, i.e., $\log S_{n,c}$ is convex. \end{conjecture}

For $c\geq 0$, U. Abel, W. Gawronski and Th. Neuschel obtained a stronger result:

\begin{theorem}\label{thm:1.2}(\cite{1}) For $c\geq 0$ the function $S_{n,c}$ is completely monotonic, i.e.,
\begin{equation*}
(-1)^j S_{n,c}^{(j)} (x)>0, \quad x\geq 0, \quad j \geq 0.
\end{equation*}

Consequently, for $c\geq 0$, $S_{n,c}$ is logarithmically convex.
\end{theorem}

The following corollary can be found in \cite{12}:

\begin{corollary}\label{cor:1.3} (\cite{12}) \begin{enumerate}[i)]
\item{} Let $c\geq 0$. Then $R_{n,c}$ is increasing and concave, while $T_{n,c}'$ is completely monotonic on $[0, +\infty )$.
\item{} $T_{n,c}$ is concave for all $c\in \mathbb{R}$.
\end{enumerate}\end{corollary}

Let us remark that the complete monotonicity for the Shannon entropy $H_{n,c}$ was investigated in \cite{12}, and for other entropies in \cite{17}.

In Sections \ref{sect:2} and \ref{sect:3} we shall use \eqref{eq:1.1} in connection with the log-convexity of $S_{n,c}$, $c\geq 0$, in order to obtain upper-bounds for $S_{n,c}$, sharper than \eqref{eq:1.3}; they can be immediately converted into sharp lower-bounds for the R\'{e}nyi entropy and the Tsallis entropy.

Theorem \ref{thm:3.1} provides an upper bound for the modified Bessel function of first kind of order $0$.

Section \ref{sect:4} is devoted to the case $c<0$. In this case Conjecture~\ref{conj:1.1} was proved in~\cite{r}, so that it is again possible to obtain upper-bounds for $S_{n,c}$, sharper than~\eqref{eq:1.3}.

On the other hand (see~\cite{6}), $S_{n,-1}$ is related to the Legendre polynomials $P_n$; using results from~\cite{6} we obtain bounds for $S_{n,-1}$ and $P_n$.

Sharp bounds on other entropies can be found in \cite{2}, \cite{4}, \cite{14}, \cite{t} and the references therein.

\section{The case $c>0$}\label{sect:2}

According to Theorem \ref{thm:1.2}, $\log{S_{n,c}(x)}$, $x\in [0,+\infty )$, is a convex function, i.e.,
\begin{equation}
S_{n,c}''(x) \geq \left ( S_{n,c}'(x) \right )^2/S_{n,c}(x), \quad x \in [0,+\infty ). \label{eq:2.1}
\end{equation}

Denote $X:=x(1+cx)$, and therefore $X'=1+2cx$. Then \eqref{eq:1.1} becomes
\begin{equation}
XX'S_{n,c}''(x)+\left ( 4(n+c)X+1 \right )S_{n,c}'(x)+2nX'S_{n,c}(x)=0. \label{eq:2.2}
\end{equation}

From \eqref{eq:2.1} and \eqref{eq:2.2} we infer that
\begin{equation}
XX'\left ( \frac{S_{n,c}'}{S_{n,c}} \right )^2 + \left ( 4(n+c)X+1 \right ) \frac{S_{n,c}'}{S_{n,c}}+2nX' \leq 0. \label{eq:2.3}
\end{equation}

This implies
\begin{equation}
\frac{S_{n,c}'(x)}{S_{n,c}(x)} \leq \frac{\sqrt{1+8cX+16(n^2+c^2)X^2}-1-4(n+c)X}{2XX'}, \label{eq:2.4}
\end{equation}
and, since $S_{n,c}(0)=1$,
\begin{equation}
\log {S_{n,c}(t)} \leq \int _0 ^t \frac{\sqrt{1+8cX+16(n^2+c^2)X^2}-1-4(n+c)X}{2XX'}\rm{dx} . \label{eq:2.5}
\end{equation}

Note that $X'^2 = 1+4cX$. Now \eqref{eq:2.5} becomes
\begin{equation}
\log {S_{n,c}(t)} \leq \int _0 ^T \frac{\sqrt{1+8cX+16(n^2+c^2)X^2}-1-4(n+c)X}{2X(1+4cX)}\rm{dX} . \label{eq:2.6}
\end{equation}
where $T:=t+ct^2$, $t\geq 0$.

Moreover, denoting $\rho := \sqrt{n^2+c^2}$ and $R:=\sqrt{16\rho ^2 T^2 + 8cT+1}$, we have

\begin{theorem}\label{thm:2.1} The following inequalities hold in the case $c>0$:
\begin{equation*}
S_{n,c}^2(t) \leq \frac{2}{1+4cT+R}\left ( \frac{1}{R+4nT} \right )^{n/c} \left ( \frac{\rho R+4\rho ^2 T+c}{\rho + c} \right )^{\rho/c} \leq
\end{equation*}
\begin{equation}
\leq \frac{1}{1+4cT}\left ( 1+4(n+c)T \right )^{-n/c} \left ( 1+8\sqrt{n^2+c^2}T \right )^{\sqrt{n^2+c^2}/c}.\label{eq:2.7}
\end{equation}

Consequently,
\begin{equation}
S_{n,c}(t) = \mathcal{O}\left ( t^{\frac{\sqrt{n^2+c^2}-n-c}{c}} \right ), \quad t\to \infty .\label{eq:2.8}
\end{equation}
\end{theorem}
\subsection*{Proof}The first inequality in \eqref{eq:2.7} follows from \eqref{eq:2.6} by a straightforward calculation. In order to get the second one it suffices to use the inequalities $1+4cT \leq R\leq 1+4\rho T$.

\begin{remark}\label{rem:2.2} The inequality \eqref{eq:2.1} is stronger than \eqref{eq:1.2}; therefore, the bound for $S_{n,c}$ given in \eqref{eq:2.7} is sharper than the bound given in \eqref{eq:1.3}. In particular, \eqref{eq:1.3} yields
\begin{equation*}
S_{n,c}(t) = \mathcal{O} \left ( t^{-\frac{n}{n+c}} \right ), \quad t\to \infty ,
\end{equation*}
and comparing with \eqref{eq:2.8} we see that
\begin{equation*}
\frac{\sqrt{n^2+c^2}-n-c}{c} < - \frac{n}{n+c}.
\end{equation*}\end{remark}

\section{The case $c=0$.\label{sect:3}}

The relations \eqref{eq:2.1} - \eqref{eq:2.6} are still valid with obvious simplifications induced by $c=0$. In particular, \eqref{eq:2.6} reduces to
\begin{equation*}
\log{S_{n,0}(t)} \leq \int _0 ^t \frac{\sqrt{1+16n^2x^2}-1-4nx}{2x}\rm{dx},
\end{equation*}
and this yields
\begin{equation}
S_{n,0}^2(t) \leq \frac{2\exp{\left ( \sqrt{1+16n^2t^2} -1-4nt\right )}}{1+\sqrt{1+16n^2t^2}} , \quad t\geq 0. \label{eq:3.1}
\end{equation}

This bound for $S_{n,0}$ is sharper than the bound furnished by \eqref{eq:1.3} with $c=0$.

By using \eqref{eq:3.1} we get also

\begin{theorem}\label{thm:3.1} Let $I_0(t)$, $t\geq 0$, be the modified Bessel function of first kind of order 0. Then
\begin{equation}
I_0^2(t) \leq \frac{2\exp{\left ( \sqrt{1+4t^2} -1 \right )}}{\sqrt{1+4t^2}+1}, \quad t\geq 0.\label{eq:3.2}
\end{equation}
\end{theorem}
\subsection*{Proof} According to \cite[(12)]{11},
\begin{equation}
I_0(t) = e^t S_{n,0}\left ( \frac{t}{2n} \right ). \label{eq:3.3}
\end{equation}
Now \eqref{eq:3.2} is a consequence of \eqref{eq:3.1} and \eqref{eq:3.3}.

\section{The case $c<0$.\label{sect:4}}

As mentioned in the Introduction, in this case Conjecture~\ref{conj:1.1} was proved in~\cite{r}. Consequently, with the same notation and the same proof as in Theorem~\ref{thm:2.1}, we get
\begin{theorem}\label{thm:4.1}
The following inequality holds for all $c<0$ and $t \in \left [ 0, -\frac{1}{c} \right ]$:
\begin{equation*}
S_{n,c}^2(t) \leq \frac{2}{1+4cT+R}\left ( \frac{1}{R+4nT} \right )^{n/c} \left ( \frac{\rho R+4\rho ^2 T+c}{\rho + c} \right )^{\rho/c}.
\end{equation*}
\end{theorem}

Since the log-convexity of $S_{n,c}$ implies the convexity, the above inequality is sharper than~\eqref{eq:1.3}. Remember that if $c<0$, then $n=-cl$ for some $l\in \mathbb{N}$. It follows that
\begin{equation*}
S_{n,c}(t)=S_{l,-1}(-ct), \quad t \in  \left [ 0, -\frac{1}{c} \right ].
\end{equation*}

Consequently, in what follows we shall investigate only the function $S_{n,-1}(x)$ with $n\in \mathbb{N}$ and $x\in [0,1]$.

G. Nikolov proved in \cite[Theorem 3]{6} that the Legendre polynomials $P_n(t)$ satisfy the inequalities
\begin{equation}
\frac{n(n+1)}{2t+(n-1)\sqrt{t^2-1}} \leq \frac{P_n'(t)}{P_n(t)} \leq \frac{n^2(2n+1)}{(n+1)t+(2n^2-1)\sqrt{t^2-1}}, \quad t\geq 1. \label{eq:4.1}
\end{equation}

Let \begin{equation*}
X:=x(1-x), \quad t=\frac{2x^2-2x+1}{1-2x} = \frac{1-2X}{X'}, \quad x\in \left [ 0,\frac{1}{2} \right ). \end{equation*}

Then $t\geq 1$ and (see \cite[(2.9)]{9}, \cite[Section 4]{11})
\begin{equation}
\frac{P_n'(t)}{P_n(t)} = \frac{nX'}{2X} + \frac{1-4X}{4X} \frac{S_{n,-1}'(x)}{S_{n,-1}(x)}. \label{eq:4.2}
\end{equation}

From \eqref{eq:4.1} and \eqref{eq:4.2} we obtain
\begin{equation}
-\frac{2nX'}{1+(n-3)X} \leq \frac{S_{n,-1}'(x)}{S_{n,-1}(x)} \leq - \frac{2n(n+1)X'}{n+1+(4n^2-2n-4)X}, \quad x \in \left [ 0,\frac{1}{2} \right ]. \label{eq:4.3}
\end{equation}

Let $t \in  \left [ 0, \frac{1}{2} \right ]$. By integrating in~\eqref{eq:4.3} with respect to $x\in [0,t]$ it follows that
\begin{equation}
\left ( 1+(n-3)T\right ) ^{-\frac{2n}{n-3}} \leq S_{n,-1}(t)\leq \left ( 1+\frac{4n^2-2n-4}{n+1}T \right )^{-\frac{n(n+1)}{2n^2-n-2}},\label{eq:4.4}
\end{equation}
where $T=t(1-t)$ and for $n=3$ the left-hand side is $e^{-6T}$. Since $S_{n,-1}(1-t) = S_{n,-1}(t)$, \eqref{eq:4.4} is valid for $t\in [0,1]$.

\begin{remark}\rm\label{rem:4.2}
For $c=-1$, \eqref{eq:1.3} is a consequence of the inequality
\begin{equation}
\frac{S'_{n,-1}(x)}{S_{n,-1}(x)}\leq -\frac{2nX'}{1+4(n-1)X}, \quad x\in \left [0, \frac{1}{2} \right ]. \label{eq:4.5}
\end{equation}
\end{remark}

Comparing~\eqref{eq:4.5} with~\eqref{eq:4.3}, we get
\begin{equation*}
-\frac{2n(n+1)X'}{n+1+(4n^2-2n-4)X} \leq - \frac{2nX'}{1+4(n-1)X}, \quad x \in \left [0, \frac{1}{2} \right ],
\end{equation*}
and so the second inequality \eqref{eq:4.4} is sharper than \eqref{eq:1.3} with $c=-1$.

\begin{remark}\rm\label{rem:4.3}
According to \cite[(29)]{11},
\begin{equation*}
S_{n,-1}(t) = \frac{1}{\pi} \int _0 ^1 \left (x+(1-x)(1-2t)^2 \right )^n \frac{dx}{\sqrt{x(1-x)}}, \quad t \in [0,1].
\end{equation*}
\end{remark}

It follows that
\begin{equation*}
S_{n,-1}(t) \geq \frac{2}{\pi} \int _0 ^1 \left (x+(1-x)(1-2t)^2 \right )^n dx,
\end{equation*}
which leads to
\begin{equation*}
\frac{1-(1-4T)^{n+1}}{2\pi (n+1)T} \leq S_{n,-1}(t), \quad t \in [0,1].
\end{equation*}

This inequality is comparable with the first inequality~\eqref{eq:4.4}.

The following results ca be found also in~\cite{9}.

Consider the inequality
\begin{equation}
\frac{P'_n(t)}{P_n(t)} \leq \frac{2n^2}{t+(2n-1)\sqrt{t^2-1}}, \quad t\geq 1, \label{eq:4.6}
\end{equation}
which was established in~\cite[Theorem 2]{6}. As remarked in~\cite{6}, \eqref{eq:4.1} is stronger than~\eqref{eq:4.6}.

From~\eqref{eq:4.6} we get by integration
\begin{equation*}
P_n(t) \leq (t+\sqrt{t^2-1})^{\frac{n(2n-1)}{2(n-1)}} \left(t+(2n-1)\sqrt{t^2-1}\right)^{-\frac{n}{2(n-1)}}, \quad t\geq 1, n\geq 2.
\end{equation*}

The stronger inequality
\begin{equation*}
P_n(t)\leq (t+\sqrt{t^2-1})^{\frac{n(2n^2-1)}{2n^2-n-2}}\left ( t + \frac{2n^2-1}{n+1} \sqrt{t^2-1} \right ) ^{-\frac{n(n+1)}{2n^2-n-2}}, \quad t\geq 1
\end{equation*}
is a consequence of~\eqref{eq:4.1}.

\section{Concluding remarks and further work}

The index of coincidence $S_{n,c}$ is intimately related with the Renyi
entropy $R_{n,c}$, Tsallis entropy $T_{n,c}$ and Legendre polynomial $P_n$. We
established new bounds for $S_{n,c}$ and, consequently, for $R_{n,c}$, $T_{n,c}$ and
$P_n$. Certain convexity properties of $S_{n,c}$ were instrumental in our
proofs. In fact, $S_{n,c}$ has also other useful convexity properties. For
example, for each integer $j$ in $[1,n]$,  $S_{n,-1}$ is $(2j-1)$-strongly convex
with modulus
\begin{equation*}
4^{j-n} {2j \choose j} {2n-2j \choose n-j}
\end{equation*}
(see the pertinent definition in \cite{a}), and for each $j\geq1$, $S_{n,0}$ is
approximately $(2j-1)$-concave with modulus
\begin{equation*}
n^{2j} {4j \choose 2j}\frac{1}{(2j)!}
\end{equation*}
(see the definition in \cite{b}).

On the other hand, according to (1.1), $S_{n,c}$ is a Heun function. By
comparing two different expressions of this Heun function it is possible
to derive combinatorial identities generalizing some classical ones from \cite{c}. Sample results are
\begin{equation*}
\sum _{j=k}^n{j \choose k}{2j \choose j}{2n-2j \choose n-j} = 4^{n-k}{n \choose k}{2k \choose k}, \quad 0\leq k \leq n,
\end{equation*}

\begin{equation*}
\sum _{i=0}^{n-j}\left ( -\frac{1}{4} \right )^i {n-j \choose i} {2i+2j \choose i+j} = 4^{j-n}{2j \choose j}{2n-2j \choose n-j}{n \choose j}^{-1}, \quad 0\leq j \leq n.
\end{equation*}

All these investigations will be presented in forthcoming papers.

\subsection*{Acknowledgements}

GM is partially supported by a grant of the Romanian Ministry of National Education
and Scientific Research, RDI Programme for Space Technology and Advanced Research -
STAR, project number 513, 118/14.11.2016.

\end{document}